\documentclass[12pt]{amsart}
\usepackage{amsfonts,amssymb,latexsym,amsmath, amsxtra}
\usepackage[all]{xy}
\usepackage[dvips]{graphics}

\usepackage[OT2,T1]{fontenc}
\DeclareSymbolFont{cyrletters}{OT2}{wncyr}{m}{n}
\DeclareMathSymbol{\Sha}{\mathalpha}{cyrletters}{"58}

\theoremstyle{plain}
\newtheorem{theorem}{Theorem}[section]

\newtheorem{lemma}[theorem]{Lemma}
\newtheorem*{conjecture}{Conjecture}
\newtheorem{proposition}[theorem]{Proposition}
\theoremstyle{definition}
\newtheorem{definition}[theorem]{Definition}
\theoremstyle{remark}
\newtheorem*{remark}{Remark}

\numberwithin{equation}{section}

\newcommand{\Z}{\mathbb Z}

\newcommand{\Q}{{\mathbb Q}}

\newcommand{\F}{{\mathbb F}}

\def\rank{\text{rank}}

\def\s2{S^{(2)}}
\def\sphi{S^{(\phi)}}
\def\shat{S^{(\hat{\phi})}}

\def\({\left(}
\def\){\right)}

\newcommand{\ol}[1]{\overline{{#1}}}

\newcommand{\mymod}[1]{\,\,({\rm mod}\,\,{#1})}
\newcommand{\logs}[1]{\frac{X (\log\log X)^{#1}}{\log X}}

\newcommand{\norm}[1]{\left|#1\right|}

\begin{document}

\title[Selmer Groups via Odd Graphs]{2-Selmer Groups and the
Birch-Swinnerton-Dyer Conjecture for the Congruent Number Curve}

\author{Robert C. Rhoades}
\address{Department of Mathematics, University of Wisconsin, Madison, WI 53706}
\email{rhoades@math.wisc.edu}

\thanks{Research of the author supported by an
NSF Graduate Research Fellowship and a National Physical Science
Consortium Graduate Research Fellowship sponsored by the NSA}

\date{\today}
\thispagestyle{empty} \vspace{.5cm}
\begin{abstract}
We take an approach toward counting the number of $n$ for which the
curves $E_n: y^2 = x^3 - n^2 x$ have 2-Selmer groups of a given
size. This question was also discussed in a pair of papers by Roger
Heath-Brown \cite{hb1, hb2}.  We discuss the connection between
computing the size of these Selmer groups and verifying cases of the
Birch and Swinnerton-Dyer Conjecture. The key ingredient for the
asymptotic formulae is the ``independence'' of the Legendre symbol
evaluated at the prime divisors of an integer with exactly $k$ prime
factors.
\end{abstract}

\maketitle 


\section{Introduction}  A  problem dating back to the tenth century,
is to determine which positive integers $n$ are  the areas of right
triangles with rational side lengths.  Such integers are called {\it
congruent numbers}. The problem of determining whether or not $n$ is
a congruent number is related to computing the rank of the elliptic
curve
\begin{equation}E_n: y^2 = x^3 - n^2x.\end{equation}  It
is well known \cite{k}
that a positive square-free integer $n$
is congruent if and only if the rank of $E_n(\Q)$, say $r(n)$, is
non-zero.
This criterion has led to many infinite families of congruent
numbers. For example, see \cite{f}, if $p$ and $q$ are distinct
primes, then we have:
\begin{itemize}
\item Heegner: $2p$ is a congruent number when $p\equiv 3
\mymod{8}$


\item Monsky: $2pq$ is congruent whenever $p \equiv 1
\mymod{8}$, $q \equiv 3, 7 \mymod{8}$ and $\(\frac{p}{q}\) = -1$.
\end{itemize}
Similarly, there are many
results that yield infinite families of non-congruent numbers.
For example, if $p,q,r, p_j$ are distinct primes, then we have:
\begin{itemize}
\item Lagrange: $pqr$ is non-congruent when $p,q \equiv 1\mymod{8}$, $r\equiv
3 \mymod{8}$, and $\(\frac{p}{q}\)= \(\frac{p}{r}\)=-1$, \cite{f},
\item Iskra: $p_1\cdots p_\ell$ is non-congruent when $p_j \equiv
3\mymod{8}$ for all $j$ and $\(\frac{p_j}{p_k}\) = -1$ for all $j
<k$, \cite{i}.
\end{itemize}
These examples show how the rank of $E_n$ is intrinsically related to
the quadratic relationships between the prime factors of $n$.
We will see this in our analysis as well.

We are interested in the size of three different Selmer groups,
which we refer to as the $2$-Selmer groups.  Let $[2]$ be the
multiplication by $2$ map, and let $\phi$ and $\hat{\phi}$ denote
the degree $2$ isogenies such that $\phi \hat{\phi} =[2]$.  More
precisely, the $2$-dual curve of $E_n$ is $E_n': y^2 =x^3 +4n^2x$
and $\phi:E_n \to E_n'$ is defined by $\phi(x,y) = (y^2/x^2,
y(n^2-x^2)/x^2)$. The $2$-Selmer groups are the Selmer groups
associated to these maps, namely, $\s2(E_n), \sphi(E_n),$ and
$\shat(E_n)$. It is a fact that $\norm{\s2(E_n)} = 2^{2+s(n)}$,
$\norm{\sphi(E_n)} = 2^{s^\phi(n)}$, and
$\norm{\shat(E_n)}=2^{2+s^{\hat{\phi}}(n)}$ for non-negative
integers $s(n), s^\phi(n)$, and $s^{\hat\phi}(n)$.
The fundamental inequality that relates $s(n)$
and $r(n)$ is
\begin{equation}r(n) \le s(n).\end{equation}
Thus information about $s(n)$ allows us to gather information about
$r(n)$.  The method of computing $\s2(E_n)$ to gain information
about the group of rational points on $E_n$ is referred to as ``full
2-descent''.  Since the rank of $\s2(E_n)$ is more accessible than
$r(n)$, it has attracted great interest in recent years (for example
see \cite{fj, f, fx, hb1, hb2, JO}).

A second approach via descent to estimating $r(n)$ is the method of
descent via isogeny.  In this approach, the sizes of $\sphi(E_n)$
and $\shat(E_n)$ are used to approximate $r(n)$.  In this approach
the fundamental inequality that relates $\sphi(E_n)$ and
$\shat(E_n)$ to $r(n)$ is
\begin{equation}
r(n) \le s^{\phi}(n)  + s^{\hat{\phi}}(n).
\end{equation}

We count the number of square-free integers up to $X$ that have
2-Selmer group of a given size. This problem was taken up and
answered precisely by Heath-Brown in a pair of papers \cite{hb1,
hb2}. Let $\lambda : = \prod_{j=1}^\infty \(1 + 2^{-j}\)^{-1}$, and
for $r = 0, 1, 2, \cdots$ let $d_r := \lambda \frac{2^r}{\prod_{1\le
j \le r} (2^j-1)}$ and let $$S(X, h) : = \{1\le n\le X: n \equiv h
\mymod{8}, n \text{ square-free}\}.$$ If $h=1$ or 3 and $r$ is even,
or if $h =5$ or 7 and $r$ is odd, then he showed that
\begin{equation}\label{eqn:heathbrown}
\norm{\{n\in S(X, h): s(n) = r\}} \sim d_r
\norm{S(X,h)}.\end{equation}

We use different techniques and obtain related results.  Heath-Brown
remarks on page 336 of \cite{hb2} that one should consider the rate
of convergence to the
 limiting distribution
in (\ref{eqn:heathbrown}) as depending on the number of prime
factors of the number $n$. For this reason we consider the problem
of determining an asymptotic for $n \in S(X,h,k)$ with 2-Selmer
groups a given size, where $S(X,h,k)$ is the set of all square-free
integers less than $X$ congruent to $h$ modulo 8 which have exactly
$k$ prime factors.

In contrast to the theorems of \cite{hb1, hb2} that deal with full
2-descent, our theorems are concerned with the descent via isogeny.
\begin{theorem}\label{thm:our HeathBrown}
In the notation above we have
\begin{align*}
\norm{\{n\in S(X, 3,k): s^{\phi}(n) = s^{\hat{\phi}}(n) =0\} } =&\,
c_3(k)(+o(1)) \norm{S(X,3,k)}\\ =& c_3(k)(1+o(1))
X\cdot \frac{(\log\log X)^{k-1}}{4 (k-1)!\log X},\\
\norm{\{n\in S(X, 2,k): s^{\phi}(n) = s^{\hat{\phi}}(n) =0\}} =&
c_2(k) (1+o(1))
\norm{S(X,2,k)}\\
=& c_2(k) (1+o(1))X\cdot \frac{(\log\log X)^{k-2}}{2 (k-2)!\log
X},\end{align*} where $c_3(k) := \frac{k}{2^{k-1}}q(k)$, $c_2(k)
:=\frac{2^{k-1} -1}{2^{2k-2}}q(k),$ and $q(k) :=
\prod_{j=1}^{\lfloor \frac{k}{2}\rfloor}
\(1-\(\frac{1}{2}\)^{2j-1}\).$
\end{theorem}


We also give a similar theorem for the Selmer groups
$\sphi(E_n)$ and $\shat(E_n)$.  Let $\omega(n)$ be the number of
prime factors of $n$.
\begin{theorem}\label{thm:HB_otherSelmer}
Let $R(X, k):= \{n < X: n \text{ square-free, } n \equiv 5 \mymod{8}
\text{ with all prime factors } \equiv 1 \mymod{4}, \omega(n) = k\}$
and $K(X,k) := \{n <X: \text{all prime factors } \equiv 1 \mymod{8},
\omega(n) =k, n \text{ square-free}\}$.  If $0\le r \le k-1$, then
\begin{align*} \norm{\{n \in R(X, k): s^{\phi}(n) =
r+1 \}}
=&\, q(k, r)(1+o(1)) \norm{R(X, k)},\\
\norm{\{n \in R(X, k): s^{\hat{\phi}}(n) = r+2\}} =&\, q(k,
r)(1+o(1)) \norm{R(X, k)},\\ \norm{\{n \in K(X, k): s^\phi(n) =
r+2\}}
=&\, q(k, r)(1+o(1)) \norm{K(X, k)},\\
\norm{\{n \in K(X, k): s^{\hat{\phi}}(n) = r+2\}} =&\, q(k,
r)(1+o(1)) \norm{K(X, k)},\end{align*} where $$q(k,s) := d(k-1,s)
\cdot 2^{\binom{k-s}{2}- \binom{k}{2}} \prod_{j=1}^{\lfloor
\frac{k-s}{2} \rfloor} \(1 - \(\frac{1}{2}\)^{2j-1}\) \text{ and }
d(m,s) = \prod_{i=0}^{r-1} \frac{2^m - 2^i}{2^s-2^i}.$$
Additionally, $\norm{R(X,k)}= (1+o(1))\frac{1}{2^{k+1} (k-1)!}X\cdot
\frac{(\log\log X)^{k-1}}{\log X}$ and $\norm{K(X,k)} = (1+o(1))
\frac{1}{4^k(k-1)!} X \cdot \frac{(\log\log X)^{k-1}}{\log X}$.
\end{theorem}



Conjecturally, the rank of an elliptic curve is
related to analytic behavior of the $L$-function associated to
$E_n$.
Let $L(E_n,s)$ be the $L$-function associated to the elliptic curve $E_n$.
Then we have the following conjecture of Birch and
Swinnerton-Dyer.
\begin{conjecture}[Birch and Swinnerton-Dyer for $E_n$]
If $E_n:y^2 = x^3 - n^2 x$, then $\text{ord}_{s=1} L(E_n, s) = r(n).$
Further if $r(n)=0$, then
\begin{equation}\label{eqn:bsd}
L(E_n, 1)/\Omega_n = 2^{\ell(n)}\norm{\Sha(E_n)},\end{equation} where
the constant $\Omega_n$ is given by \begin{equation}\label{eqn:period}
\Omega_n : = \frac{1}{\sqrt{n}}
\int_1^\infty (x^3-x)^{-1/2}dx\end{equation} and
$\ell(n)$ is a non-negative integer.  Also, $\Sha(E_n)$ is the
Tate-Shafarevich group of $E_n$ over $\Q$.
\end{conjecture}
A special case of a famous theorem of Rubin \cite{r1, r2},
implies that since $E_n$ has complex multiplication by $\Z[i]$,
if $L(E_n, 1) \ne 0$, then the
group $\Sha(E_n/\Q)$ is finite and the odd parts of both sides of
equation (\ref{eqn:bsd}) are equal. Thus, Rubin's result
proves the Birch
and Swinnerton-Dyer conjecture for $E_n$ with $r(n) =0$ up to a
power of $2$.
For many cases we are able to compute the power of 2 appearing in both
sides of equation (\ref{eqn:bsd}) and show that they are equal. Hence,
we
conclude the truth of the Birch and Swinnerton-Dyer conjecture for
the elliptic curve $E_n$, for many values of $n$.


An advantage of the approach used here, that is restricting the
number of prime factors, is that we can analyze the analytic
properties of the $L$-functions at the same time we study the
arithmetic properties of $E_n$.  This is a result of the fact that
the same combinatorial conditions used to analyze the size of the
2-Selmer groups appear in the analysis of the $2$-power in the
$L$-value. In fact, Zhao in a series of papers \cite{z1, z2, z3, z4}
described these conditions.  Combining our results with Zhao's work
and the work of Feng and Xiong \cite{fx}, gives the following:
\begin{theorem}\label{thm:BSD1}
Let $q(k)$ be as in Theorem \ref{thm:our HeathBrown}.
\begin{enumerate}
\item Let $B_3(X, k)$ be the set of all $n<X$ with $\omega(n) =k$,
$n\equiv  3 \mymod{8}$, where $n$ has exactly one prime factor
congruent to 3 modulo 8, all other prime factors are 1 modulo 8. For
any $k$, the Birch and Swinnerton-Dyer Conjecture is true for all
$n\in B_3(\infty, k)$ with $s^\phi(n)= s^{\hat{\phi}}(n) =0$.
Moreover, we have that
$$\norm{\{n\in B_3(X, k): s^\phi(n)= s^{\hat{\phi}}(n)
= 0 \}} = q(k)(1+o(1))\norm{B_3(X,k)}.$$  Additionally,
$\norm{B_3(X,k)} = (1+o(1))\frac{k}{4^k (k-1)!} X\cdot
\frac{(\log\log X)^{k-1}}{\log X}.$

\item Let $B_2(X, k)$ be the set of all $n<X$ with $n\equiv  2 \mymod{8}$,
$\omega(n/2) =k$, $n$ has all odd prime factors congruent to 1
modulo 4. For any $k$, the Birch and Swinnerton-Dyer Conjecture is
true for all $n\in B_2(\infty, k)$ with $s^\phi(n)=
s^{\hat{\phi}}(n) =0$.  Moreover, we have that $$\norm{\{n\in B_2(X,
k): s^\phi(n)= s^{\hat{\phi}}(n) = 0 \}} = \frac{2^k-1}{2^k}q(k)
(1+o(1))\norm{B_2(X,k)}.$$ Additionally, $\norm{B_2(X,k)} = (1+o(1))
\frac{1}{2^{k+1}(k-1)!} X \cdot \frac{(\log \log X)^{k-1}}{\log X}.$
\end{enumerate}
\end{theorem}

\begin{remark}
This theorem gives information about the number of twists of $L(E_1,
s)$  which have $L(E_n, 1) \ne 0$. Ono and Skinner \cite{o, os} have
given much more general results which establish lower bounds for the
number of twists of an $L$-function which have non-vanishing of the
central value.  Similar to this result, their results are based on
showing the ``oddness'' of the algebraic part of the $L$-value.
\end{remark}

Theorem \ref{thm:BSD1} verifies BSD for all curves $E_n$ with
trivial 2-Selmer groups $\sphi(E_n)$ and $\shat(E_n)$ such that the
prime factors of $n$ are subject to some congruence conditions.  It
is possible to use the work of Zhao \cite{z4} and Li and Tian
\cite{lt} to verify the full BSD conjecture for an infinite class of
curves $E_n$ whose Tate-Shafarevich group has non-trivial 2-part.

\begin{theorem}\label{thm:BSD2}
Let $D_1(X, k)$ be the set of all $n<X$ with $\omega(n) =k$,
$n\equiv  1 \mymod{8}$, $n$ has all prime factors congruent to 1
modulo 8.  Then for any $k$, the Birch and Swinnerton-Dyer
Conjecture is true for all $n\in D_1(\infty, k)$ with $s^\phi(n)
=2$, $s^{\hat{\phi}}(n) =0$ and $\Sha(E_n)[2]=\Z/2\Z\times \Z/2\Z$.
Also if the conditions on the Selmer groups and Tate-Shafarevich
group are satisfied then $r(n) =0$.  Finally,
\begin{align*}&\norm{\{n\in D_1(X, k): s^\phi(n) =2, s^{\hat{\phi}}(n)
=0 ,\Sha(E_n)[2] =\Z/2\Z\times \Z/2\Z  \}}\\ & \,\,\ge q(k)
\frac{1}{2}(1+o(1))\norm{D_1(X,k)} = \frac{q(k)}{2^{2k+1}
(k-1)!}(1+o(1)) X \cdot \frac{(\log \log X)^{k-1}}{\log X}.
\end{align*}
\end{theorem}

Briefly our approach to showing the vanishing of the $\phi$ and
$\hat{\phi}$-Selmer groups for the curve $E_n$ follows from
congruence conditions on the prime factors of $n$ along with the
equidistribution of the Legendre symbol.  Feng \cite{f} introduced
the language of ``odd graphs'' to encode the necessary information
about the Legendre symbol.  In Section \ref{sec:ODD GRAPHS} we
introduce his language and some related results that we will need in
our analysis.  Section \ref{SEC:square-free} puts congruence
conditions on a classical theorem of Landau which gives an
asymptotic for the number of integers less than $X$ with a fixed
number of prime factors as $X$ tends to infinity.  The
equidistribution of the Legendre symbol that is needed essentially
follows from work of Cremona and Odoni \cite{co} and is recalled in
Section \ref{SEC:orthogLegendre}.  In Section \ref{sec:HEATH-BROWN}
we combine the results of Sections \ref{sec:ODD
GRAPHS}-\ref{SEC:orthogLegendre} to give the proofs of Theorems
\ref{thm:our HeathBrown} and
\ref{thm:HB_otherSelmer}. 
Finally, in Section \ref{sec:BSD} we use Zhao's  results \cite{z1,
z2, z3, z4} to verify some cases of the Birch and Swinnerton-Dyer
Conjecture.

%

\section*{Acknowledgments}
I am indebted to Roger Heath-Brown and an anonymous referee for
comments on earlier versions of this paper. I thank Ken Ono for many
useful discussions and comments. I would also like to thank Neil
Calkin and Kevin James for originally introducing this problem to me
during an REU at Clemson University during the summer of 2004.
\section{Counting Selmer Groups Via Odd Graphs}
\label{sec:ODD GRAPHS}

The theory of ``odd graphs'', initiated by Feng \cite{f}, has been
used in many places to count Selmer groups, see \cite{fj, f, fx}.
This section describes two theorems of Feng and Xiong \cite{fx} that
gives necessary and sufficient conditions for the triviality of the
$\phi$ and $\hat{\phi}$-Selmer groups. We will also describe some
results of Faulkner and James \cite{fj} which we will use to prove
Theorem \ref{thm:HB_otherSelmer}. We now describe the graphs we will
be interested in.

Throughout this section, unless otherwise stated, suppose that $n$
is an odd square-free integer with $n = p_1\dots p_t>0$.  Define the
directed graphs $G(n)$, $G(-n)$, $G'(n)$ by
\begin{equation}
V(G(n)) := \{p_1, \ldots, p_t\} \hspace{.05in} \text{ and }
\hspace{.05in} E(G(n)) := \bigg\{\overrightarrow{p_ip_j}:
\(\frac{p_i}{p_j}\) = -1, 1 \le i\ne j \le t \bigg\},
\end{equation}
\begin{align}
V(G(-n)) := &\{-1,p_1, \ldots, p_t\} \hspace{.05in} \text{ and } \\
 E(G(-n)) := &\bigg\{\overrightarrow{p_ip_j}: \(\frac{p_i}{p_j}\) =
-1, 1 \le i\ne j \le t, p_i\neq 3\pmod{4}\bigg\} \\& \cup
\bigg\{\overrightarrow{(-1)r}: r\in \{p_1, \cdots, p_r\}, r\equiv
\pm 3 \mymod{8}\bigg\},\nonumber
\end{align}
\begin{align}
V(G'(n)) := &\{2,p_1, \ldots, p_t\} \hspace{.05in} \text{ and } \\
 E(G'(n)) := &\bigg\{\overrightarrow{p_ip_j}: \(\frac{p_i}{p_j}\) =
-1, 1 \le i\ne j \le t, p_i\neq 3\pmod{4}\bigg\} \\& \cup
\bigg\{\overrightarrow{r2}: r\in \{p_1, \cdots, p_r\}, r\equiv \pm 5
\mymod{8}\bigg\},\nonumber
\end{align}
where $V(\cdot)$ and $E(\cdot)$ stand for the vertex set and edge
sets of the graph.

\begin{definition}
Suppose that $G$ is a graph with vertex set $V$ and edge set $E$. A
{\em partition} of $G$ is a pair $(S,T)$ of sets such that $S \cap T
= \emptyset$ and $S\cup T = V$. A partition $(S, T)$ is {\it even}
provided that all $v\in S$ have an even number of edges directed
from $v$ to vertices in $T$ and all $v\in T$ have an even number of
edges directed from $v$ to vertices in $S$.
\end{definition}
In particular, the partitions $(G, \emptyset)$ and $(\emptyset, G)$
are always even partitions. We call these the trivial partitions.
Let $e(G)$ be the number of even partitions of the graph $G$.
\begin{definition}
A graph $G$ is called {\em even} provided that it admits a
nontrivial even partition.  A graph $G$ is said to be {\em odd}
provided that its only even partitions are trivial.
\end{definition}
We recall two theorems of
\cite{fx} that we will use to obtain Theorem \ref{thm:our HeathBrown}.

\begin{theorem}[Theorem 2.4 of \cite{fx}]\label{thm:2.4 of fx}
Suppose that $n \equiv \pm 3\mymod{8}$.  Then $\sphi(E_n) = \{1\}$
and $\shat(E_n) = \{\pm 1,\pm n\}$ if and only if the following
three conditions are satisfied:\begin{enumerate}

\item $n \equiv 3\mymod{8}$
\item $n = p_1 \dots p_t$, $p_1 \equiv 3\mymod{4}$ and $p_j \equiv 1
\mymod{4}$ for  $(2\le j \le t)$.
\item  $G(n)$ is an odd graph. \end{enumerate}
\end{theorem}


\begin{theorem}[Theorem 2.6 of \cite{fx}]\label{thm:2.6 of fx}
Suppose that $2 \parallel n$ then $\sphi(E_n) = \{1\}$ and
$\shat(E_n) = \{\pm 1, \pm n\}$ if and only if $G'(n/2)$ is odd.
Furthermore, if $\sphi(E_n) = \{1\}$ and $\shat(E_n) = \{\pm 1 , \pm
n\}$ then all odd primes dividing $n$ are 1 modulo 4 and there is at
least one that is 5 modulo 8.
\end{theorem}

We now state the results of Faulkner and James that we will use. The
following combines special cases of Theorems 1.4 and 1.5 of
\cite{fj}.
\begin{theorem}\label{thm:fj}
With the notation from above, we have \begin{enumerate}
\item If $n\equiv 5 \mymod{8}$,
and $n$ has all primes congruent to 1 modulo 4, then
$$\norm{\sphi(E_n)} =  e(G(n )) ,$$ and
$$\norm{\shat(E_n)} = 2\cdot e(G(n)).$$

\item If $n\equiv 1 \mymod{8}$ and all the prime factors of $n$ are
congruent to 1 modulo 8, then $$\norm{\sphi(E_n)} = 2\cdot
e(G(n)),$$ and
$$\norm{\shat(E_n)} = e(G(-n)).$$
\end{enumerate}
\end{theorem}

\subsection{Approach For Asymptotics}
The approach for using Theorems \ref{thm:2.4 of fx} and \ref{thm:2.6
of fx} to prove Theorem \ref{thm:our HeathBrown} is the following:
For any $n$ that satisfies the first two conditions of Theorem
\ref{thm:2.4 of fx}, by quadratic reciprocity $G(n)$ is an
undirected graph. Furthermore, by Dirichlet's theorem on primes in
arithmetic progressions and induction on the number of prime factors
of $n$, we see that for any undirected graph $G$ there exist
infinitely many $n$ such that $G(n) = G$.

Given that there are infinitely many $n$ such that each undirected
graph $G$ appears as $G(n)$, one might hope that selecting $n$ (that
satisfies the first two conditions of Theorem \ref{thm:2.4 of fx})
at random would result in selecting a random graph $G(n)$. To make
this more precise we fix an integer $k$. Say $\{n_1, n_2, n_3,
\dots, n_M\}$ is the set of all appropriate integers less than some
$X$ with $k$ prime factors. Then we might hope that if we look at
the list of graphs $(G(n_1), G(n_2), \dots, G(n_M))$ then each
undirected graph on $k$ vertices appears in the list with equal
proportion.  If this is true, then the proportion of $n\in \{n_1,
n_2, \dots, n_M\}$ that have $S^{(\phi)}(n)$ and
$S^{(\hat{\phi})}(n)$ trivial should be the same as the probability
that a random undirected graph on $k$ vertices is odd.


\subsection{Probability of Odd Graphs}
Before moving on, we recall the results of \cite{us} about the
probability that a graph on $k$ vertices is odd and some results
connecting the number of odd partitions of a graph and the rank of
an associated matrix over
$\F_2$. 

If $G=(V,E)$ is a graph with vertices $v_1, \cdots, v_k$, then
define the adjacency matrix $A(G)$ of a graph $G$ by $A(G) =
(a_{ij})_{1 \le i, j \le k}$ where for $i\ne j$, $a_{ij} = 1$ if
$\overrightarrow{v_iv_j}\in E(G)$ and 0 otherwise and $a_{ii} =0$.
Let $d_i = \sum_{j=1}^k a_{ij}\mymod{2}$.  The Laplace matrix of $G$
is defined by $L(G) = \text{diag}(d_1, \cdots, d_k) + A(G)$.
\begin{lemma}[Lemma 2.2 \cite{fx}]\label{lem:graphToMatrix}
Let $G = (V, E)$ be a directed graph, $k = \norm{V}$ and $r=
\rank_{\F_2} L(G)$.  Then the number of even partitions of $G$ is
$2^{k-r}$.  In particular, $G$ is an odd graph if and only if
$r=k-1$.
\end{lemma}
\begin{remark}
In fact, \cite{fj} shows an explicit relationship between elements
of the kernel of the matrix $L(G)$ and elements of the Selmer groups
under consideration.
\end{remark}
\begin{theorem}[Theorem 1.6 \cite{us}]\label{thm:probGraphOdd}
Let $G$ be an undirected graph on $k$ vertices. Denote the
probability that $G$ has $2^{e+1}$ even partitions by $q(k,e)$, for
$0 \le e \le k-1$. Then
$$q(k,e) =2^{\binom{k-e}{2} - \binom{k}{2}} d(k-1,e)
\prod_{j=1}^{\lfloor \frac{k-e}{2} \rfloor}
\(1-\(\frac{1}{2}\)^{2j-1}\),
$$
where $d(m,j) := \prod_{i=0}^{j-1} \frac{2^m-2^i}{2^j - 2^i}.$
\end{theorem}

Also in \cite{us} one finds the following proposition.
\begin{proposition}\label{cor:probRankFull}
Denote the probability that a $k\times k$ matrix over $\F_2$ has
rank $k$ by $p(k)$.  Then $p(k) = q(k+1, 0) = \prod_{j=1}^{\lfloor
(k+1)/2 \rfloor} (1-2^{-2j+1})$, where $q(k+1,k)$ is as in Theorem
\ref{thm:probGraphOdd}.
\end{proposition}

We will need the following proposition as well.
\begin{proposition}\label{prop:moreOdd}
Let $A$ be a $k\times k$ symmetric matrix over $\F_2$.  Given that
the sum of the rows of $A$ is $v = (0, \cdots, 0, 1 ,\cdots, 1)^T$
for some vector with $j>0$ 1's, the probability that $A$ has rank
$k$ is $p(k-1)=q(k,0)$, independent of $j$.
\end{proposition}
\begin{proof}
If $j=1$, then apply the remark after Theorem 2.6 of
\cite{fx}.  If $j\ge 2$, then we can find a change of basis matrix
$\Lambda$ such that $\Lambda v = (0\cdots 0 \, 1)^T$.  Thus
$\Lambda^T A \Lambda$ is a symmetric matrix with the same rank as
$A$ and we are reduced to the $j=1$ case.
\end{proof}

\begin{remark}
Monsky, in an appendix to \cite{hb2}, gives a way to compute the
size of $S^{(2)}(n)$ by computing the $\F_2$ rank of larger matrix
than any of the ones we consider here.
\end{remark}


\section{Square-free Integers with factors in Specified Congruence
Classes}\label{SEC:square-free}

In this section we answer the question of how many integers $n < X$
have exactly $k$ prime factors where each lies in a specified
congruence class modulo $m$, for some $m$. The most classical result
in this direction, due to Landau, states:
\begin{equation}\label{eqn:landau}\left|\{n\le X: \Omega(n) = k
\}\right| \sim \left|\{n\le X: \omega(n) = k \} \right| \sim \frac{X
(\log\log X)^{k-1}}{(k-1)!\log X},\end{equation} where $\Omega(n)$
is the number of prime factors of $n$ counted with multiplicity and
$\omega(n)$ is the number of prime factors of $n$ counted without
multiplicity.

We fix the following notation: for a fixed positive integer $m>1$,
we let $1= r_1< \cdots < r_{\phi(m)} <m$ be the $\phi(m)$ standard
representatives for $(\Z/m\Z)^\times$.  Define $\pi_k(X;m;a_1,
\cdots, a_{\phi(m)})$ to be the number of $n\le X$ with $\omega(n) =
k$ and $n$ is square-free with exactly $a_j$ prime factors congruent
to $r_j \mymod{m}$ for $1 \le j \le \phi(m)$. We have the following
theorem.
\begin{theorem}\label{thm:square-free}
Let $k$ and $m$ be fixed positive integers with $m>1$ and let $X$ be
a positive real number.  If  $0\le a_1, a_2, \cdots, a_{\phi(m)} \le k$ are
integers such that $a_1 + a_2 + \cdots + a_{\phi(m)} = k$, then
$$\pi_k(X;m;a_1, a_2, \cdots, a_{\phi(k)}) = \(1+o(1)\)
\frac{k!}{a_1! a_2! \cdots a_{\phi(m)}!} \frac{1}{\phi(m)^k (k-1)!}
\frac{X(\log\log X)^{k-1}}{\log(X)}.$$
\end{theorem}

%

We will prove Theorem \ref{thm:square-free} by induction on the
number of non-zero $a_j$.  The base case is where just one of the
$a_j$ is non-zero.  As short hand we write $\pi_k(X;m;b) =
\pi_k(X;m;a_1, \cdots , a_{\phi(m)} )$ where there is a $j$ with
$r_j = b$ and $a_j = k$. In this case we necessarily have $a_i=0$
for all $i\ne j$.
\begin{proposition}\label{prop:square-free ONE}
If $m$ is a positive integer greater than 1 and $b$ is a positive
integer relatively prime to $m$, then
$$\pi_k(X;m;b) = \(1+o(1)\)\frac{1}{\phi(m)^k(k-1)!}
\frac{X (\log\log X)^{k-1}}{\log(X)}.$$
\end{proposition}

This result follows from a trivial modification of the proof of
Landau's result.  For a proof of Landau's result see
\cite{nathanson1}.
To complete the induction we will need the following two lemmas.
\begin{lemma}\label{lem:recipricol sum}
Let $k$ be a positive integer, and $C$ a positive real number.  Let
$S$ be a set of positive integers such that $$S(X):=\norm{ \{n\le X:
n\in S\} } = \(1+o(1)\) \frac{C}{(k-1)!} \frac{X (\log\log
X)^{k-1}}{\log
(X) },$$ as $X\to \infty$.  
We have that $$\int_{2}^{X^{1/2}} \frac{S(t)}{t^2 \log(X/t)} dt =\(1+o(1)\)
\frac{C}{k!} \frac{(\log\log X)^{k}}{\log (X)}.$$
\end{lemma}
\begin{proof}
For any $\epsilon >0$, there exists an  $N$ such that
for all $X >N$, $$\norm{S(X)-  \frac{C}{(k-1)!}\logs{k-1}} <
\epsilon \frac{C}{(k-1)!} \logs{k-1}.$$  Therefore it follows that
\begin{align*}
\int_{2}^{X^{1/2}} \frac{S(t)}{t^2 \log(X/t)} dt =& \int_2^N
 \frac{S(t)}{t^2 \log(X/t)} dt + \int_{N}^{X^{1/2}}
\frac{S(t)}{t^2 \log(X/t)} dt \\
=& O\(\frac{\log N}{\log X}\) + \int_{N}^{X^{1/2}} \frac{S(t)}{t^2
\log(X/t)} dt.
\end{align*}
Where we estimate the first integral by using the fact that $S$ is a
set of positive integers and so $S(t) \le t$.

Now we turn to estimating $\int_{N}^{X^{1/2}} \frac{S(t)}{t^2
\log(X/t)} dt$.  In the range of integration we know that we may
replace $S(t)$ by its asymptotic formula and introduce a small
error.  Precisely we have
$$\norm{\int_{N}^{X^{1/2}}
\frac{S(t)}{t^2 \log(X/t)} dt - \frac{C}{(k-1)!}\int_{N}^{X^{1/2}}
\frac{t(\log\log t)^{k-1}}{t^2\log(t)  \log(X/t)} dt} \le \epsilon
\frac{C}{(k-1)!} \int_{N}^{X^{1/2}} \frac{t (\log\log
t)^{k-1}}{t^2\log(t)  \log(X/t)} dt.$$ We have
\begin{align*}
\int_{N}^{X^{1/2}} \frac{(\log\log t)^{k-1}}{t \log(t) \log(X/t)} dt
=& \frac{1}{\log(X)}\int_{N}^{X^{1/2}} \frac{(\log\log t)^{k-1}}{t
\log(t)}\(1+O\(\frac{\log(t)}{\log(X)}\)\) dt \\
=& \frac{(\log\log X)^k}{k \log(X)} + O\(\frac{1}{(\log X)^2}
\int_N^{X^{1/2}} \frac{(\log\log t)^{k-1}}{t} dt\) \\
=& \frac{(\log\log X)^k}{k \log(X)} + O\(\frac{(\log\log
X)^{k-1}}{\log X}\).
\end{align*}
Thus we may conclude that
$$\int_{2}^{X^{1/2}} \frac{S(t)}{t^2 \log(X/t)} dt =
\frac{C(\log\log X)^k}{k! \log(X)}(1+o(1)) + O\(\frac{\log N}{\log
X}\).$$  Taking $N$ to be of size $\log\log X$ is sufficient to
yield the conclusion of the lemma.
\end{proof}

\begin{lemma}\label{lem:prod_disjoint_sets}
Let $k_1, k_2$ be positive integers and $C_1$ and $C_2$ be positive
numbers.  Let $S_1$ and $S_2$ be two sets of positive integers such
that for each $j$
$$S_j(X):=\norm{ \{n\le X: n\in
S_j\} } = \(1+o(1)\) \frac{C_j}{(k_j-1)!} \frac{X (\log\log
X)^{k_j-1}}{\log (X) },$$ as $X\to \infty$.  We have that
$$S_{1,2}(X):= \norm{\{ (n_1, n_2) \in S_1\times S_2 : n_1 n_2 \le X\}}
= \(1+o(1)\) \frac{C}{(k_1+k_2-1)!} \frac{X(\log\log X)^{k_1+k_2 -
1}}{\log X},$$ where $C = \frac{C_1C_2
(k_1+k_2-1)!}{(k_1-1)!(k_2-1)}\(\frac{1}{k_1} + \frac{1}{k_2}\)$.
\end{lemma}

\begin{proof}We abuse notation and refer to $S_j(X)$ as the set
of elements in $S_j$ which are less than or equal to $X$, as well as
the size of the set of elements of $S_j$ up to $X$.  The use of the
symbol will be clear from the context.

Begin with the following inclusion-exclusion-like identity
$$ S_{1,2}(X) = \sum_{t \in S_2(X^{1/2})} S_1(X/t) + \sum_{t\in
S_1(X^{1/2})} S_2(X/t) -  S_1(X^{1/2}) S_2(X^{1/2}).$$ This identity
follows from the fact that the first two sums will count everything
in the set $S_{1,2}(X)$, while over counting precisely those
elements which equal $n_1n_2$ where $n_1 \in S_1$ and $n_2 \in S_2$
and both $n_1$ and $n_2$ are smaller than $X^{1/2}$.

By assumption we have $$S_1(X^{1/2})S_2(X^{1/2}) = O\( \frac{X
(\log\log X)^{k_1 + k_2 -2}}{(\log X)^2}\),$$  and this is well within
our expected error. 
We now estimate the first sum.

Since $t \le X^{1/2}$ we know that $X/t \ge X^{1/2}$, so we may
apply our asymptotic to obtain
$$\sum_{t \in S_2(X^{1/2})} S_1(X/t) = (1+o(1)) \frac{C_1}{(k_1-1)!}X
\sum_{t \in S_2(X^{1/2})}\frac{(\log\log X/t)^{k_1
-1}}{t\log(X/t)}.$$ Using the fact that $\log\log(X/t) = \log\log X
+ \log\(1-\frac{\log t}{\log X}\) = \log\log X + O(1)$ for $t\in
[1,X^{1/2}]$, we obtain
\begin{equation}\label{eqn:inclusion_sum}\sum_{t \in S_2(X^{1/2})}
S_1(X/t) = (1+o(1)) \frac{C_1}{(k_1-1)!}X (\log\log X)^{k_1 -1}
\sum_{t \in S_2(X^{1/2})}\frac{1}{t\log(X/t)}.\end{equation} We have
\begin{align*} \sum_{t \in S_2(X^{1/2})}\frac{1}{t\log(X/t)} =&
\int_2^{X^{1/2}}
\frac{1}{t\log(X/t)}d S_2(t) \\
=& \frac{S_2(X^{1/2})}{X^{1/2} \log(X^{1/2})} - \int_2^{X^{1/2}}
\frac{S_2(t)}{t^2 (\log(X/t))^2} dt + \int_2^{X^{1/2}}
\frac{S_2(t)}{t^2 \log(X/t)} dt\\
=& O\(\frac{(\log\log X)^{k_2-1}}{(\log X)^2}\) -O\( \frac{(\log\log
X)^{k_2}}{(\log X)^2}\) + \(1+o(1)\) \frac{C_2}{k_2!}
\frac{(\log\log X)^{k_2}}{\log (X)} \\
=&\(1+o(1)\) \frac{C_2}{k_2!} \frac{(\log\log X)^{k_2}}{\log (X)},
\end{align*}
where we use the fact that
$$ \int_2^{X^{1/2}}
\frac{S_2(t)}{t^2 (\log(X/t))^2} dt \ll \frac{1}{\log X}
\int_2^{X^{1/2}} \frac{S_2(t)}{t^2 \log(X/t)} dt$$ and Lemma
\ref{lem:recipricol sum} to estimate the first integral.

Combining this with equation (\ref{eqn:inclusion_sum}) we obtain
\begin{equation*}\sum_{t \in S_2(X^{1/2})} S_1(X/t) = (1+o(1))
\frac{C_1C_2}{(k_1-1)!k_2!} \logs{k_1+k_2-1}.
\end{equation*}
The exact same argument for the sum $\sum_{t \in S_1(X^{1/2})}
S_2(X/t)$, shows that
$$S_{1,2}(X) = (1+o(1))
\(\frac{C_1C_2}{(k_1-1)!k_2!}  + \frac{C_1C_2}{k_1! (k_2
-1)!}\)\logs{k_1+k_2-1}.$$
\end{proof}

With this lemma in hand we prove Theorem \ref{thm:square-free}.
\begin{proof}[Proof of Theorem \ref{thm:square-free}]
Using Proposition \ref{prop:square-free ONE} as the base case,  the
result now follows by induction with a straightforward application
of Lemma \ref{lem:prod_disjoint_sets}.
\end{proof}

\section{Independence of Legendre Symbols}\label{SEC:orthogLegendre}
The main theorem of this section is to prove that the Legendre
symbols are independent; this allows us to conclude that the graphs
discussed in Section \ref{sec:ODD GRAPHS} are asymptotically uniformly
distributed.
This theorem is a simple extension of the results in Section 3 of
\cite{co}. As a result we do not include the proof here.  Instead we
only remark that to obtain this theorem one would follow the
argument of \cite{co} but would need to add the additional
constraint on the $L$-functions considered to take into account the
additional congruence conditions on the prime factors $p_j$.
Specifically, we have the following theorem:
\begin{theorem}\label{thm:LegendreInd}
Let $k\ge 2$ be a positive integer.  Fix $\epsilon_{ij} \in \{-1,
1\}$ and $\delta_j \in \{1, 3, 5, 7\}$ for $1\le j \le k$ and $1\le
i< j\le k$.  For ease of notation let $\delta = (\delta_1, \cdots,
\delta_k)$.  Let $C_k(X, \delta)$ be the set of $k$-tuples $(p_1,
\cdots, p_k)$ of primes with $2<p_1 < p_2 < \cdots < p_k \le X$,
$p_1 \cdots p_k \le X$, $p_j \equiv \delta_j \mymod{8}$.  Then the
number of elements of $C_k(X, \delta)$ with $\(\frac{p_i}{p_j}\) =
\epsilon_{ij}$ for $i <j$, is
$$2^{-\binom{k}{2}}(1+o(1)) \norm{C_k(X,\delta)}.$$
\end{theorem}
\section{Selmer Group Asymptotics}\label{sec:HEATH-BROWN}
In this section we give the proofs of Theorems \ref{thm:our
HeathBrown} and \ref{thm:HB_otherSelmer}.  The strategy for the
proofs was explained in Section \ref{sec:ODD GRAPHS}.  Here we
quickly give the proofs, which amount to combining the results from
the previous sections.

\begin{proof}[Proof of Theorem \ref{thm:our HeathBrown}]
Begin with the case $n=p_1p_2\cdots p_k \equiv 3 \mymod{8}$. Then by
Theorem \ref{thm:2.4 of fx} and Lemma \ref{lem:graphToMatrix}, we
know that $s^\phi(n)= s^{\hat{\phi}}(n) =0$ if and only if  $n=p_1
p_2 \cdots p_k$ with $p_1 \equiv 3 \mymod{4}$, $p_j \equiv 1
\mymod{4}$, and $G(n)$ has exactly 2 even partitions (after suitable
renaming the $p_j$). By Theorem \ref{thm:square-free}, we know that
the number of $n\equiv 3 \mymod{8}$ with the necessary congruence
conditions on the prime factors is $\frac{k}{2^{k-1}}
(1+o(1))\norm{S(X,3,k)}$. By Theorems \ref{thm:probGraphOdd} and
\ref{thm:LegendreInd}, we know that of all the $n$ with the
necessary congruence conditions on the prime factors the proportion
of them with $G(n)$ odd is $q(k,0)(1+o(1)) = q(k)(1+o(1))$.   The
result follows.

The proof for $n\equiv 2 \mymod{8}$ is similar.  However instead of
appealing to Theorem \ref{thm:probGraphOdd} we use Proposition
\ref{prop:moreOdd}.   From Theorem \ref{thm:2.6 of fx}, we know that
if $s^\phi(n)= s^{\hat{\phi}}(n) =0$ for some $n\equiv 2\mymod{8}$,
then there are no primes congruent to 3 modulo 4 that divide $n$.
Say $n = 2p_1 \cdots p_{k-1}$ has $k-1$ odd prime factors. Hence
$L(G'(n))$ is a $k \times k$ matrix with
$$L(G'(n)) = \left(\begin{array}{cc} A & v \\
0 \cdots 0 & 0  \end{array} \right),$$ where $A$ is a symmetric
$(k-1)\times (k-1)$ matrix determined by $\(\frac{p_i}{p_j}\)$ and
the vector $v$ is a matrix with the same number of 1's as primes
congruent to 5 modulo 8, which by Theorem \ref{thm:2.6 of fx} is
necessarily larger than zero.  Now we know that the sum of the rows
of $L(G'(n))$ is 0.  So we know, by Lemma \ref{lem:graphToMatrix}
that the graph $G'(n)$ is odd if and only if the matrix $A$ has full
rank.  We apply Proposition \ref{prop:moreOdd} to see the
probability that $G'(n)$ is odd is $q(k-1,0)=q(k-1)$. Finally, we
apply Theorem \ref{thm:LegendreInd} to justify that each possible
$k\times k$ matrix appears with equal probability.
\end{proof}

\begin{proof}[Proof of Theorem \ref{thm:HB_otherSelmer}]
Our starting point is Theorem \ref{thm:fj}. This proof is similar to
the proof of Theorem \ref{thm:our HeathBrown}, however it is
important to realize that because the conditions on the prime
factors of $n$ the graphs are all undirected or equivalently all the
matrices $L(G(n))$ or $L(G(-n))$ are symmetric matrices.  Therefore
we may apply Theorem \ref{thm:probGraphOdd} to determine the
probability that the matrix has a given rank.  We apply Theorem
\ref{thm:LegendreInd} to see that it is appropriate to treat the
graphs appearing for such $n$ as random undirected graphs.
\end{proof}


\section{Verifying BSD}\label{sec:BSD}
In this section we prove Theorems \ref{thm:BSD1} and \ref{thm:BSD2}.
The proof of the first of these theorems amounts to combining
results from the work of Feng and Xiong, \cite{fx}, the work of Zhao
\cite{z1, z3, z4}, and the results from this paper.  The second
theorem uses some work of Zhao \cite{z2} and work of Li and Tian
\cite{lt}.

Because the work of Zhao is important we will state one of his three
theorems which we will employ. Let $L(\ol{\psi}_{n^2}, s)$
denote the Hecke $L$-function corresponding to the dual of
$\psi_{n^2}$ which is the Gr\"{o}ssencharacter of $\Q[i]$ attached
to $E_n$.  
See \cite{z1} for more details. Let $\Omega_n$ be as in equation
(\ref{eqn:period}).

\begin{theorem}[Theorem 2 of \cite{z1}]
Suppose $n = p_1 \cdots p_m$ with $p_1 \equiv 3 \mymod{8}$ and $p_j
\equiv 1 \mymod{8}$ for all $j>1$.  The power of 2 in
$L(\ol{\psi}_{n^2}, 1)/\Omega_n$ is  greater than or equal to
$2m-1$ with equality if and only if
$G(n)$ is odd.
\end{theorem}
Recall that for $n$ with the prime factorization of this theorem we
know from Theorem \ref{thm:2.4 of fx} that $s(n) =0$ and thus
$\Sha(E_n)$ is odd when $G(n)$ is odd.  Since the condition for the
lowest power of 2 is the same here as it is for $\Sha(E_n)$ to be
odd we are able to verify the Birch and Swinnerton-Dyer Conjecture.
Indeed, Zhao gives:
\begin{proposition}[Proposition 3 of
\cite{z1}]\label{prop:3of20} Suppose $n \equiv 3 \mymod{8}$, $n$ has
one prime factor congruent to 3 modulo 8 and all others congruent to
1 modulo 8. If $G(n)$ is odd, then the Birch and Swinnerton-Dyer
Conjecture is true.
\end{proposition}

\begin{proof}[Proof of Theorem \ref{thm:BSD1}]
Theorem \ref{thm:2.4 of fx} shows that $n\in B_3(X,k)$, $s^\phi(n)=
s^{\hat{\phi}}(n)=0$ if and only if the graph $G(n)$ is odd.
Proposition \ref{prop:3of20} shows that for such $n$ the Birch and
Swinnerton-Dyer conjecture is true. Finally, applying Theorem
\ref{thm:probGraphOdd} and Theorem \ref{thm:LegendreInd}, we obtain
the asymptotic, as in the proof of Theorem \ref{thm:our HeathBrown}.
To prove the second case of this theorem, we use Theorem
\ref{thm:2.6 of fx} and Corollary 3 of \cite{z3}.
\end{proof}

The same proof using the main theorem of \cite{z4} with Theorem 2.5
of \cite{fx} gives the following proposition.
\begin{proposition}
\item Let $B_1(X, k)$ be the set of all $n<X$ with $\omega(n) =k$,
$n\equiv  1 \mymod{8}$, $n$ has exactly two prime factors congruent
to 3 modulo 8, all other prime factors are 1 modulo 8.  Then for any
$k$, the Birch and Swinnerton-Dyer Conjecture is true for all $n\in
B_1(\infty, k)$ with $s^\phi(n)= s^{\hat{\phi}}(n) =0$.
\end{proposition}
Since the combinatorics of this case are a bit messier we do not
bother to present the asymptotics.   Before we can prove the final
result of this paper, we must give one more definition and one more
lemma.

\begin{definition}
Let $p\equiv 1 \mymod{8}$ be prime.  Then set $\delta(p) = 1$ if we
have one of the following:
\begin{itemize}
\item $p\equiv 1 \mymod{16}$ and $\(\frac{2}{p}\)_4 =-1$
\item $p \equiv 9 \mymod{16}$ and $\(\frac{2}{p}\)_4 = 1$,
\end{itemize}
and set $\delta (p) = 0$ otherwise.   Here $\(\frac{2}{p}\)_4$ is
the quartic character.  For an integer $n= p_1 \cdots p_k$ with each
$p_j \equiv 1 \mymod{8}$, set $\delta(n) = \sum_{i=1}^k \delta(p_i)
\mymod{2}$.
\end{definition}
The following result is important for Theorem \ref{thm:BSD2}
\begin{proposition}[Li and Tian \cite{lt}]\label{prop:LiTianMore}
Let $n\in D_1(X,k)$, where $D_1$ is as in Theorem \ref{thm:BSD2}. We
have $G(n)$ is an odd graph and $\delta(n) =1$ then $s^\phi(n) =2$,
$s^{\hat{\phi}}(n) =0$ , and $\Sha(E_n)[2] = \Z/2\Z \times \Z/2\Z$ .
\end{proposition}
The converse direction should also be true, but for brevity we do
not give the proof. Proving the converse direction would result in
an equality in Theorem \ref{thm:BSD2} rather than the lower bound
given.

\begin{theorem}[Theorem 1 of \cite{z2}]\label{thm:zhoa_2}
Suppose $n\in D_1(X,k)$.  Then the power of 2 in
$L(\ol{\psi}_{n^2},1)/\Omega_n \ge 2k$ and there is equality if and
only if $\delta(n)$ is odd and $G(n)$ is an odd graph.
\end{theorem}

\begin{proof}[Proof of Theorem \ref{thm:BSD2}]
For $n\in D_1(X,k)$ we proceed as in the proof of Theorem
\ref{thm:BSD1}.  Zhao \cite{z4} proved that the power of 2 in the
$L$-value is as small as possible if and only if $G(n)$ is odd and
$\delta(n)$ is odd.  The previous proposition proves that the power
of 2 in $\Sha(E_n)$ is as small as possible (in this case 2) if and
only if $\delta(n)$ is odd and $G(n)$ is odd.  Now we obtain the
asymptotic by applying Theorem \ref{thm:LegendreInd} combined with
Theorem \ref{thm:probGraphOdd} and noting that for $\norm{\{n\in
D_1(X, k): \delta(n) = 1 \}} = (1/2+o(1)) \norm{D_1(X,k)}$.
Technically we would need a version of Theorem
\ref{thm:LegendreInd}, which has $\delta_j \in \{1,9\}$ and we
consider the $p_j$ modulo 16 instead of 8.  But the proof goes
through the same as the case we handle.
\end{proof}

\end{document}